\documentclass{article}
\usepackage{amsmath}
\usepackage{amsfonts}
\usepackage{amssymb}

\newtheorem{theorem}{Theorem}

\newtheorem{corollary}[theorem]{Corollary}

\newtheorem{definition}[theorem]{Definition}

\newtheorem{lemma}[theorem]{Lemma}

\newtheorem{remark}[theorem]{Remark}

\newcommand{\R}{{\mathbb{R}}}

\begin{document}

\title{Variational representations for $N$-cyclically monotone vector fields}
\author{Alfred Galichon\thanks{%
This research has received funding from the European Research Council under
the European Union's Seventh Framework Programme (FP7/2007-2013) / ERC grant
agreement n%
${{}^\circ}$%
313699. Support from FiME, Laboratoire de Finance des March\'{e}s de
l'Energie (www.fime-lab.org) is gratefully acknowledged.} \\
\textit{\small Economics Department}\\
\textit{\small Sciences Po Paris }\\
\textit{{\small 28 rue des Saints-P\`{e}res,}} \textit{\small 75007 Paris
France }\\
\textit{\small alfred.galichon@sciences-po.fr} \and Nassif Ghoussoub\thanks{%
Partially supported by a grant from the Natural Sciences and Engineering
Research Council of Canada.} \hspace{2mm} \\
\textit{\small Department of Mathematics}\\
\textit{\small University of British Columbia}\\
\textit{\small Vancouver BC Canada V6T 1Z2}\\
\textit{\small nassif@math.ubc.ca}\vspace{1mm}}
\date{July 10, 2012; Last revised October 2, 2013}
\maketitle

\begin{abstract}
Given a convex bounded domain $\Omega $ in ${{\mathbb{R}}}^{d}$ and an
integer $N\geq 2$, we associate to any \textit{jointly $N$-monotone} $(N-1)$-tuplet $(u_1, u_2,..., u_{N-1})$ of vector fields from $%
\Omega$ into $\mathbb{R}^{d}$, a Hamiltonian $H$ on ${\mathbb{R}}^{d} \times
{\mathbb{R}}^{d} ... \times {\mathbb{R}}^{d}$, that is concave in the first
variable, jointly convex in the last $(N-1)$ variables such that for almost all $%
x\in \Omega$,
\begin{equation*}
\hbox{$(u_1(x), u_2(x),..., u_{N-1}(x))= \nabla_{2,...,N} H(x,x,...,x)$.}
\end{equation*}%
Moreover, $H$ is $N$-sub-antisymmetric, meaning that $\sum%
\limits_{i=0}^{N-1}H(\sigma ^{i}(\mathbf{x}))\leq 0$ for all $\mathbf{x}%
=(x_{1},...,x_{N})\in \Omega ^{N}$, $\sigma $ being the cyclic permutation
on ${\mathbb{R}}^{d}$ defined by $\sigma
(x_{1},x_2,...,x_{N})=(x_{2},x_{3},...,x_{N},x_{1})$. Furthermore, $H$ is $N$%
-antisymmetric in a sense to be defined below.
This can be seen as an extension of a theorem of E. Krauss, which associates
to any monotone operator, a concave-convex antisymmetric saddle function. We
also give various variational characterizations of vector fields that are
almost everywhere $N$-monotone, showing that they are dual to the class of
measure preserving $N$-involutions on $\Omega$. .
\end{abstract}

\section{Introduction}

Given a domain $\Omega $ in ${{\mathbb{R}}}^{d}$, recall that a
single-valued map $u$ from $\Omega $ to ${{\mathbb{R}}}^{d}$ is said to be
\textit{$N$-cyclically monotone} if for every cycle $%
x_{1},...,x_{N},x_{N+1}=x_{1}$ of points in $\Omega $, one has%
\begin{equation}
\sum_{i=1}^{N}\left\langle u\left( x_{i}\right) ,x_{i}-x_{i+1}\right\rangle
\geq 0.
\end{equation}%
A classical theorem of Rockafellar \cite{Ph} states that  a map $u$ from $%
\Omega$ to ${\mathbb{R}}^d$ is \textit{$N$-cyclically monotone for every $%
N\geq 2$} if and only if
\begin{equation}
\hbox{$u (x)\in \partial \phi (x)$ for all $x\in \Omega$,}
\end{equation}
where $\phi:{\mathbb{R}}^d \to {\mathbb{R}}$ is a convex function. On the
other hand, a result of E. Krauss \cite{Kra} yields that $u$ is a monotone
map, i.e., a $2$-cyclically monotone map,  if and only if
\begin{equation}
\hbox{$u(x)\in \partial_2H(x, x)$ for all $x\in \Omega$,}
\end{equation}
where $H$ is a concave-convex antisymmetric Hamiltonian on ${\mathbb{R}}%
^d\times {\mathbb{R}}^d$, and $\partial_2H$ is the subdifferential of $H$ as
a convex function in the second variable.

In this paper, we extend the result of Krauss to the class of $N$-cyclically
monotone vector fields, where $N\geq 3$. We shall give a representation for
a family of $(N-1)$ vector fields, which may or may not be individually $N$%
-cyclically monotone. Here is the needed concept.

\begin{definition}
Let $u_1,..., u_{N-1}$ be bounded vector fields from a domain $\Omega
\subset {\mathbb{R}}^d$ into ${\mathbb{R}}^d$. We shall say that the $(N-1)$%
-tuple $(u_1, u_2,..., u_{N-1})$ is \textit{jointly $N$-monotone}, if for
every cycle $x_{1},...,x_{2N-1}$ of points in $\Omega $ such that $%
x_{N+i}=x_{i}$ for $1\leq i\leq N-1$, one has%
\begin{equation}
\sum_{i=1}^{N} \sum_{\ell=1}^{N-1} \langle u_{l}(x_i), x_i-x_{l+i} \rangle
\geq 0.
\end{equation}
\end{definition}

\textbf{Examples of jointly $N$-monotone families of vector fields:}

\begin{itemize}

\item It is clear that $(u,0,0,...,0)$ is jointly $N$-monotone if and only
if $u$ is $N$-monotone.

\item More generally, if each $u_{\ell }$ is $N$-monotone, then the family $%
(u_{1},u_{2},...,u_{N-1})$ is jointly $N$-monotone. Actually, one only needs
that for $1\leq \ell \leq N-1$, the vector field $u_{\ell }$ be \textit{$%
(N,\ell )$-monotone}, in the following sense: for every cycle $%
x_{1},...,x_{N+\ell }$ of points in $\Omega $ such that $x_{N+i}=x_{i}$ for $%
1\leq i\leq \ell $, we have
\begin{equation}
\sum_{i=1}^{N}\left\langle u_{\ell }\left( x_{i}\right) ,x_{i}-x_{\ell
+i}\right\rangle \geq 0.
\end{equation}%
This notion is sometimes weaker than $N$-monotonicity since
if $\ell $ divides $N$, then it suffices for $u$ to be $\frac{N}{\ell }$%
-monotone in order to be an $(N,\ell )$-monotone vector field. For example,
if $u_{1}$ and $u_{3}$ are $4$-monotone operators and $u_{2}$ is $2$%
-monotone, then the triplet $(u_{1},u_{2},u_{3})$ is jointly $4$-monotone.

\item Another example is when $(u_{1},u_{2},u_{3})$ are vector fields such
that $u_{2}$ is $2$-monotone and
\begin{equation*}
\hbox{$\langle u_1(x)-u_3(y), x-y\rangle \geq 0$ for every $x, y \in \R^d$.}
\end{equation*}%
In this case, the triplet $(u_{1},u_{2},u_{3})$ is jointly $4$-monotone. In
particular, if $u_{1}$ and $u_{2}$ are both $2$-monotone, then the triplet $%
(u_{1},u_{2},u_{1})$ is jointly $4$-monotone.

\item More generally, it is easy to show that $(u,u,...,u)$ is jointly $N$-monotone if and only if $u$ is $2$-cyclically monotone.
\end{itemize}

\noindent In the sequel, we shall denote by $\sigma $ the cyclic permutation
on ${\mathbb{R}}^{d}\times ...\times {\mathbb{R}}^{d}$, defined by
\begin{equation*}
\sigma (x_{1},x_2,...,x_{N-1},x_{N})=(x_{2},x_{3},...,x_{N},x_{1}),
\end{equation*}%
and consider the family of continuous \emph{$N$-antisymmetric}
Hamiltonians on $\Omega ^{N}$, that is
\begin{equation}\label{exact}
\mathcal{H}_{N}(\Omega )=\{H\in C(\Omega
^{N});\,\sum\limits_{i=0}^{N-1}H(\sigma ^{i}(x_{1},...,x_{N}))=0\}.
\end{equation}%
We say that $H$ is \emph{$N$-sub-antisymmetric on $\Omega$}
if
\begin{equation}\label{subanti}
\hbox{$\sum\limits_{i=0}^{N-1}H(\sigma ^{i}(x_{1},...,x_{N}))\leq 0 $ on $%
\Omega^N$.} 
\end{equation}
We shall also say that a function $F$ of two variables is \emph{$N$-cyclically sub-antisymmetric on $\Omega$}, if
\begin{equation}
\hbox{$F(x,x)=0$ and $\sum\limits_{i=1}^{N} F(x_i, x_{i+1})\leq 0$ for all cyclic families $x_1,..., x_N, x_{N+1}=x_1$ in $\Omega$.}
\end{equation}
Note that if a function $H(x_1,..., x_N)$ $N$-sub-antisymmetric and if it only depends on the first two variables, then the function $F(x_1, x_2):=H(x_1, x_2,..., x_N)$
is $N$-cyclically sub-antisymmetric.

We associate to any function $H$ on $\Omega ^{N}$, the following functional
on $\Omega \times ({\mathbb{R}}^{d})^{N-1}$,
\begin{equation}
L_{H}(x,p_{1},...,p_{N-1})=\sup \left\{ \sum_{i=1}^{N-1}\langle
p_{i},y_{i}\rangle -H(x,y_{1},...,y_{N-1});y_{i}\in \Omega \right\} .
\label{legendre*}
\end{equation}%
Note that if $\Omega $ is convex and if $H$ is convex in the last $(N-1)$
variables, then $L_{H}$ is nothing but the Legendre transform of $\tilde{H}$
with respect to the last $(N-1)$ variables, where $\tilde{H}$ is the
extension of $H$ over $({\mathbb{R}}^d)^N$, defined as: $\tilde{H}=H$ on $%
\Omega ^{N} $ and $\tilde{H}=+\infty $ outside of $\Omega ^{N}$. Since $%
H(x,...,x)=0$ for any $H\in {\mathcal{H}}_{N}(\Omega)$, then for any such $H$%
, we have for $x\in \Omega$ and $p_1,..., p_{N-1} \in {\mathbb{R}}^d$,
\begin{equation}
L_{H}(x,p_{1},...,p_{N-1})\geq \sum_{i=1}^{N-1}\langle x,p_{i}\rangle .
\end{equation}
To formulate variational principles for such vector fields, we shall
consider the class of $\sigma$-invariant probability measures on $\Omega^N$,
which are those $\pi \in {\mathcal{P}}(\Omega^N)$ such that for all $h\in
L^1(\Omega^N, d\pi)$, we have
\begin{equation}
\int_{\Omega^N} h(x_1,..., x_N) d\pi=\int_{\Omega^N} h(\sigma (x_1,...,
x_N)) d\pi.
\end{equation}
We denote
\begin{equation}
\hbox{${\mathcal P}_{\rm sym}(\Omega^N)=\{\pi \in {\mathcal P}(\Omega^N)$; $\pi$ $\sigma$-invariant  probability
on $\Omega^N\}$.}
\end{equation}
For a given probability measure $\mu$ on $\Omega$, we also consider the
class
\begin{equation}
\hbox{${\mathcal P}^{\mu}_{{\rm sym}}(\Omega^N)=\{\pi \in {\mathcal P}_{\rm
sym}(\Omega^N)$;  ${\rm proj}_1{\pi}=\mu$\}},
\end{equation}
i.e., the set of all $\pi \in {\mathcal{P}}_{\mathrm{sym}}(\Omega^N)$ with a
given first marginal $\mu$, meaning that
\begin{equation}
\hbox{$\int_{\Omega^N} f(x_1)\, d\pi (x_1,..., x_N) =\int_{\Omega} f(x_1)\, d\mu
(x_1)$ for every $f\in L^1(\Omega, \mu)$.}
\end{equation}
Consider now the set ${\mathcal{S}}(\Omega,\mu)$ of $\mu$-measure preserving
transformations on $\Omega$, which can be identified with a closed subset of
the sphere of $L^2(\Omega, {{\mathbb{R}}}^d)$. We shall also consider the
subset of ${\mathcal{S}}(\Omega, \mu)$ consisting of $N$-involutions, that
is
\begin{equation*}
\hbox{${\mathcal S}_N(\Omega, \mu)=\huge\{S\in {\mathcal S}(\Omega, \mu);\,  S^N=I$
$\mu$ a.e.\}.}
\end{equation*}

\section{Monotone vector fields and $N$-antisymmetric Hamiltonians}

In this section, we establish the following extension of a theorem of Krauss.

\begin{theorem}
\label{main.1} Let $N\geq 2$ be an integer, and consider $u_1,..., u_{N-1}$
to be bounded vector fields from a convex domain $\Omega \subset {\mathbb{R}}%
^d$ into ${\mathbb{R}}^d$.

\begin{enumerate}
\item If the $(N-1)$-tuple $(u_1,..., u_{N-1})$ is jointly $N$-monotone,
then there exists an $N$-sub-antisymmetric Hamiltonian $H$ that is zero on the diagonal of $\Omega^N$, concave
in the first variable, convex in the other $(N-1)$ variables such that
\begin{equation}  \label{rep}
\hbox{$(u_1(x),..., u_{N-1}(x))=\nabla_{2,...,N} H(x,x,...,x)$ \quad
for a.e. $x\in \Omega$.}
\end{equation}%
Moreover, $H$ is $N$-antisymmetric in the following sense
\begin{equation}
H(x_{1},x_{2},...,x_{N})+H_{2,...,N}(x_{1},x_{2},...,x_{N})=0,
\label{symmH1}
\end{equation}%
where $H_{2,...,N}$ is the concavification of the function $K(\mathbf{x}%
)=\sum\limits_{i=1}^{N-1}H(\sigma ^{i}(\mathbf{x}))$ with respect to the
last $(N-1)$ variables.

Furthermore, there exists a continuous $N$-antisymmetric Hamiltonian ${\bar H}$ on $\Omega ^{N}$, such that
\begin{equation}\label{dualrep}
\hbox{$L_{\bar H}(x, u_1(x), u_2(x),..., u_{N-1}(x))=\sum\limits_{i=1}^{N-1}\langle u_i(x), x\rangle$ for all $x\in \Omega$.}
\end{equation}

\item Conversely, if $(u_1,..., u_{N-1})$ satisfy (\ref{rep}) for some $N$%
-sub-antisymmetric Hamiltonian $H$ that is zero on the diagonal of $\Omega^N$, concave in the first variable,
convex in the other variables, then the $(N-1)$-tuple $(u_1,..., u_{N-1})$
is jointly $N$-monotone on $\Omega$.
\end{enumerate}
\end{theorem}

\begin{remark}
Note that in the case $N=2$, $K\left( \mathbf{x}\right) =H\left(
x_{2},x_{1}\right) $ is concave with respect to $x_{2}$, hence $H_{2}\left(
x_{1},x_{2}\right) =H\left( x_{2},x_{1}\right) $, and (\ref{symmH1}) becomes%
\begin{equation*}
H\left( x_{1},x_{2}\right) +H\left( x_{2},x_{1}\right) =0,
\end{equation*}%
thus $H$ is antisymmetric, recovering well-known results \cite{Kra}, \cite%
{Gh}, \cite{G-M}, \cite{Mi}.
\end{remark}

We start with the following lemma.

\begin{lemma}
\label{cute} Assume the $(N-1)$-tuple of bounded vector fields $(u_1,...,
u_{N-1})$ on $\Omega$ is jointly $N$-monotone. Let $f(x_1,...,
x_N):=\sum_{l=1}^{N-1} \langle u_{l}(x_1), x_1-x_{l+1} \rangle $ and
consider the function ${\tilde f}(x_1,...,x_n)$ to be the convexification of
$f$ with respect to the first variable, that is
\begin{equation}
{\tilde f}\left( x_1, x_2,..., x_N\right) =\inf \left\{
\sum_{k=1}^{n}\lambda _{k}f\left( x^k_{1},x_2,..., x_N\right): \, n \in {%
\mathbb{N}}, \lambda _{k}\geq 0,\sum_{k=1}^{n}\lambda
_{k}=1,\sum_{k=1}^{n}\lambda _{k}x_{1}^k=x_1\right\} .
\end{equation}
Then, ${\tilde f}$ satisfies the following properties:

\begin{enumerate}
\item $f\geq {\tilde f}$ on $\Omega^N$;

\item ${\tilde f}$ is convex in the first variable and concave with respect
to the other variables;

\item ${\tilde f}(x,x,..., x)=0$ for each $x\in \Omega$,

\item ${\tilde f}$ satisfies
\begin{equation}  \label{good}
\hbox{$\sum_{i=0}^{N-1} {\tilde f}(\sigma^i(x_1,...,x_N))\geq 0$ on
$\Omega^N$.}
\end{equation}
\end{enumerate}
\end{lemma}

\noindent\textbf{Proof:} Since the $(N-1)$-tuple $(u_1,..., u_{N-1})$ is
jointly $N$-monotone, it is easy to see that the function
\begin{equation*}
f(x_1,..., x_N):=\sum_{l=1}^{N-1} \langle u_{l}(x_1), x_1-x_{l+1} \rangle
\end{equation*}
is linear in the last $(N-1)$ variables, that $f(x,x,...,x)=0$, and that
\begin{equation}
\hbox{$\sum_{i=0}^{N-1} f(\sigma^i(x_1,...,x_N))\geq 0$ on $\Omega^N$.}
\end{equation}
It is also clear that $f\geq {\tilde f}$, that ${\tilde f}$ is convex with
respect to the first variable $x_1$, and that it is concave with respect to
the other variables $x_2,..., x_N$, since $f$ itself is concave (actually
linear) with respect to $x_2,..., x_N$. We now show that $\tilde f$
satisfies (\ref{good}).

For that, we fix $x_{1},x_{2},...,x_{N}$ in $\Omega $ and consider $%
(x_{1}^{k})_{k=1}^{n}$ in $\Omega $, and $(\lambda _{k})_{k}$ in ${\mathbb{R}%
}$ such that $\lambda _{k}\geq 0$ such that $\sum_{k=1}^{n}\lambda _{k}=1$
and $\sum_{k=1}^{n}\lambda _{k}x_{1}^{k}=x_{1}$. For each $k$, we have
\begin{equation*}
f(x_{1}^{k},x_{2},..., x_N)+ f(x_{2},..., x_N, x_{1}^k)+ ...+f(x_N, x_{1}^k,
x_{2},..., x_{N-1}) \geq 0.
\end{equation*}%
Multiplying by $\lambda_k$, summing over $k$, and using that $f$ is linear
in the last $(N-1)$-variables, we have
\begin{equation*}
\sum_{k=1}^{n}\lambda _{k}f(x_{1}^{k},x_{2},..., x_N)+ f(x_{2},..., x_N,
x_{1})+ ...+f(x_N, x_{1}, x_{2},..., x_{N-1}) \geq 0.
\end{equation*}%
By taking the infimum, we obtain
\begin{equation*}
{\tilde f}\left( x_{1},x_{2},..., x_N\right) +\sum_{i=1}^{N-1}f(\sigma^i
(x_1,x_2,..., x_N))\geq 0.
\end{equation*}%
Let now $n\in {\mathbb{N}}$, $\lambda _{k}\geq 0$, $x_{N}^{k}\in \Omega $ be
such that $\sum_{k=1}^{n}\lambda _{k}=1$ and $\sum\limits_{k=1}^{n}\lambda
_{k}x_{2}^{k}=x_{2}$. We have for every $1\leq k\leq n$,
\begin{equation*}
{\tilde f}\left( x_{1},x_{2}^k, x_3,..., x_N\right) + f\left(x_2^k,
x_{3},,..., x_1\right)+...+ f\left(x_N, x_{1},x_{2}^k, x_3,...,
x_{N-1}\right) \geq 0.
\end{equation*}%
Multiplying by $\lambda _{k}$, summing over $k$ and using that $\tilde f$ is
convex in the first variable and $f$ is linear in the last $(N-1)$%
-variables, we obtain
\begin{eqnarray*}
{\tilde f}\left( x_{1},x_{2}, x_3,..., x_N\right) +
\sum\limits_{k=1}^{n}\lambda _{k}f\left(x_2^k, x_{3},,..., x_1\right)+...+
f\left(x_N, x_{1},x_{2}, x_3,..., x_{N-1}\right) \\
\geq \sum\limits_{k=1}^{n}\lambda _{k}{\tilde f}\left( x_{1},x_{2}^k,
x_3,..., x_N\right) + \sum\limits_{k=1}^{n}\lambda _{k}f\left(x_2^k,
x_{3},,..., x_1\right)+...+ \sum\limits_{k=1}^{n}\lambda _{k} f\left(x_N,
x_{1},x^k_{2}, x_3,..., x_{N-1}\right) \geq 0.
\end{eqnarray*}%
By taking the infimum over all possible such choices, we get
\begin{equation*}
{\tilde f}\left( x_{1},x_{2}, x_3,..., x_N\right) + {\tilde f}\left(x_2,
x_{3},,..., x_1\right) +...+f\left(x_N, x_{1},x_{2}, x_3,..., x_{N-1}\right)
\geq 0.
\end{equation*}
By repeating this procedure with $x_{3},...,x_{N-1}$, we get
\begin{equation*}
\sum_{i=0}^{N-2}{\tilde f}\left(\sigma^i(x_1, x_{2},,..., x_N)\right)+
f\left(x_N, x_{1},x_{2}, x_3,..., x_{N-1}\right) \geq 0.
\end{equation*}
Finally, since
\begin{equation*}
f\left(x_N, x_{1},x_{2}, x_3,..., x_{N-1}\right) \geq -\sum_{i=0}^{N-2}{%
\tilde f}\left(\sigma^i(x_1, x_{2},,..., x_N)\right).
\end{equation*}
and since ${\tilde f}$ is concave in the last $(N-1)$ variables, we have for
fixed $x_{1},x_{2},...,x_{N-1}$, that the function
\begin{equation*}
x_{N}\rightarrow -\sum_{i=0}^{N-2}{\tilde f}\left(\sigma^i(x_1, x_{2},,...,
x_N)\right)
\end{equation*}
is a convex minorant of $x_{N}\rightarrow f\left(x_N, x_{1},x_{2}, x_3,...,
x_{N-1}\right)$. It follows that
\begin{eqnarray*}
f\left(x_N, x_{1},x_{2}, x_3,..., x_{N-1}\right)\geq {\tilde f}\left(x_N,
x_{1},x_{2}, x_3,..., x_{N-1}\right) \geq -\sum_{i=0}^{N-2}{\tilde f}%
\left(\sigma^i(x_1, x_{2},,..., x_N)\right),
\end{eqnarray*}%
which finally implies that $\sum_{i=0}^{N-1}{\tilde f} (%
\sigma^i(x_1,x_2,...,x_N))\geq 0.$

This clearly implies that ${\tilde f}(x,x,...,x)\geq 0$ for any $x\in \Omega$%
. On the other hand, since ${\tilde f}(x,x,...,x)\leq f(x,x,...,x)=0$, we
get that ${\tilde f}(x,x,...,x)=0$ for all $x\in \Omega$. \hfill $\Box$
\bigskip

\noindent \textbf{Proof of Theorem \ref{main.1}:} Assume the $(N-1)$-tuple
of vector fields $(u_1,..., u_{N-1})$ is jointly $N$-monotone on $\Omega$,
and consider the function $f(x_1,..., x_N):=\sum_{l=1}^{N-1} \langle
u_{l}(x_1), x_1-x_{l+1} \rangle $ as well as its convexification with
respect to the first variable ${\tilde f}(x_1,...,x_N)$.

By Lemma \ref{cute}, the function $\psi (x_1,...,x_N):=-{\tilde f}%
(x_1,...,x_N)$ satisfies the following properties
\begin{trivlist}

\item (i) $x_1\to \psi (x_1,...,x_N)$ is concave;
\item (ii) $(x_2,x_3,..., x_N)\to \psi (x_1,...,x_N)$ is convex;
\item (iii) $\psi (x_1,...,x_N)\geq -{f}(x_1,...,x_N)=\sum_{l=1}^{N-1} \langle u_{l}(x_1), x_{l+1}-x_1 \rangle $;
\item (iv) $\psi $ is $N$-sub-antisymmetric.
\end{trivlist}
Consider now the family $\overline{\mathcal{H}}$ of functions $H: \Omega^N\to {%
\mathbb{R}}$ such that

\begin{enumerate}
\item $H(x_1,x_2,..., x_N) \geq \sum_{l=1}^{N-1} \langle u_{l}(x_1),
x_{l+1}-x_1 \rangle $ for every $N$-tuple $(x_1,..., x_N)$ in $\Omega^N$;

\item $H$ is concave in the first variable;

\item $H$ is jointly convex in the last $(N-1)$ variables;

\item $H$ is $N$-sub-antisymmetric.
\item $H$ is zero on the diagonal of $\Omega^N$.
\end{enumerate}

Note that $\overline{\mathcal{H}}\neq \emptyset$ since $\psi$ belongs to $\overline{\mathcal{H}}$. Note that any $H$ satisfying (1) and (4) automatically satisfies (5). Indeed, by $N$-sub-antisymmetry, 
we have for all
$\mathbf{x}=(x_1,..., x_N)\in \Omega^N$,
\begin{equation}
H(\mathbf{x})\leq - \sum_{i=1}^{N-1}H(\sigma^i(\mathbf{x}))\leq
-\sum_{i=1}^{N-1} \psi (\sigma^i(\mathbf{x)}).
\end{equation}
This also yields that
\begin{equation}  \label{bounds}
\sum_{\ell=1}^{N-1}\langle u_\ell (x_1), x_{\ell +1}-x_1\rangle \leq H(%
\mathbf{x})\leq -\sum_{i=2}^{N} \sum_{\ell=1}^{N-1}\langle u_\ell(x_{i}),
x_{i}- x_{i+\ell}\rangle,
\end{equation}
where we denote $x_{i+N}:=x_i$ for $i=1,...,\ell$. This yields that $%
H(x,x,...,x)=0$ for any $x\in \Omega$.

It is also easy to see that every directed family $(H_i)_i$ in
$\overline{\mathcal{H}}$ has a supremum $H_\infty\in \overline{\mathcal{H}}$, meaning that $\overline{{
\mathcal{H}}}$ is a Zorn family, and therefore has a maximal element $H$.

Consider now the function
\begin{equation*}
\bar H (\mathbf{x})=\frac{(N-1)H(\mathbf{x})- \sum_{i=1}^{N-1}H(\sigma^i(%
\mathbf{x}))}{N},
\end{equation*}
and note that\newline

\noindent (i) $\bar H$ is $N$-antisymmetric, since
\begin{equation*}
\bar H(\mathbf{x})=\frac{1}{N}\sum_{i=1}^{N-1}[H(\mathbf{x})-H(\sigma^i (%
\mathbf{x}))],
\end{equation*}
and each $K_i(\mathbf{x}):=H(\mathbf{x})-H(\sigma^i (\mathbf{x}))$ is $N$%
-antisymmetric.\newline

\noindent (ii) $\bar H \geq H$ on $\Omega^N$, since
\begin{equation*}
N[\bar H(\mathbf{x}) - H(\mathbf{x})]=- \sum_{i=0}^{N-1}H(\sigma^i(\mathbf{x}%
)) \geq 0,
\end{equation*}
because $H$ itself is $N$-sub-antisymmetric.\newline

The maximality of $H$ would have implied that $H=\bar{H}$ is $N$%
-antisymmetric if only $\bar{H}$ was jointly convex in the last $(N-1)$%
-variables, but since this is not necessarily the case, we consider for $%
\mathbf{x}=(x_{1},x_{2},...,x_{N})$, the function
\begin{equation*}
K(x_{1},x_{2},...,x_{N})=K(\mathbf{x}):=-\sum_{i=1}^{N-1}H(\sigma ^{i}(%
\mathbf{x})),
\end{equation*}%
which is already concave in the first variable $x_{1}$. Its convexification
in the last $(N-1)$-variables, that is
\begin{equation*}
K^{2,...,N}(\mathbf{x})=\inf \left\{ \sum\limits_{i=1}^{n}\lambda _{i}{K}%
(x_{1},x_{2}^{i},...,x_{N}^{i});\,\lambda _{i}\geq
0,\sum\limits_{i=1}^{n}\lambda
_{i}(x_{2}^{i},...,x_{N}^{i},1)=(x_{2},...,x_{N},1)\right\} ,
\end{equation*}%
is still concave in the first variable, but is now convex in the last $(N-1)$
variables. Moreover,
\begin{equation}
H\leq K^{2,...,N}\leq K=-\sum_{i=1}^{N-1}H\circ \sigma ^{i}.
\end{equation}%
Indeed, $K^{2,...,N}\leq K$ from the definition of $K^{2,...,N}$, while $%
H\leq K^{2,...,N}$ because $H\leq K$ and $H$ is already convex in the last $%
(N-1)$-variables. It follows that
\begin{equation*}
H\leq \frac{(N-1)H+K^{2,...,N}}{N}\leq \frac{(N-1)H+K}{N}=\frac{%
(N-1)H-\sum\limits_{i=1}^{N-1}H\circ \sigma ^{i}}{N}=\bar{H}.
\end{equation*}%
The function $H^{\prime }=\frac{(N-1)H+K^{2,...,N}}{N}$ belongs to the
family $\overline{\mathcal{H}}$ and therefore $H=H^{\prime }$ by the maximality of $H$%
.

This finally yields that $H$ is $N$-sub-antisymmetric, that $H(x,...,x)=0$ for
all $x\in \Omega $ and that
\begin{equation*}
\hbox{$H({\bf x}) +H_{2,..., N}({\bf x})=0$ for every ${\bf x}\in \Omega^N$,}
\end{equation*}%
where $H_{2,...,N}=-K^{2,...,N}$, which for a fixed $x_1$, is nothing but
the concavification of $(x_2, ..., x_N) \rightarrow \sum_{i=1}^{N-1}H(\sigma
^{i}(x_1,x_2, ..., x_N))$.

Note now that since for any $x_1,..., x_N$ in $\Omega$,
\begin{equation}  \label{inequa}
H(x_1,x_2, ... x_N) \geq \sum_{\ell=1}^{N-1}\langle u_\ell (x_1), x_{\ell
+1}-x_1\rangle,
\end{equation}
and
\begin{equation}  \label{diago}
H(x_1,x_1,..., x_1)=0,
\end{equation}
we have
\begin{equation}  \label{conv}
H(x_1,x_2, ..., x_N) -H(x_1,..., x_1) \geq \sum_{\ell=1}^{N-1}\langle u_\ell
(x_1), x_{\ell +1}-x_1\rangle.
\end{equation}
Since $H$ is convex in the last $(N-1)$ variables, this means that for all $%
x\in \Omega$, we have
\begin{equation}  \label{rep1}
(u_1(x), u_2(x),...,u_{N-1}(x))\in \partial_{2,..., N}H(x, x,..., x).
\end{equation}
as claimed in (\ref{rep}).
Note that this also yields that
\begin{equation*}
\hbox{$ L_H(x,u_1(x),...,u_{N-1}(x)) +
H(x,x,...,x)=\sum_{\ell=1}^{N-1}\langle u_\ell(x), x\rangle$ for all
$x\in \Omega$. }
\end{equation*}
In other words, $L_H(x,u_1(x),...,u_{N-1}(x))=\sum_{\ell=1}^{N-1} \langle u_\ell(x), x\rangle$ \, for all $x\in \Omega$.
As above, consider
\begin{equation*}
\bar H (\mathbf{x})=\frac{(N-1)H(\mathbf{x})- \sum_{i=1}^{N-1}H(\sigma^i(%
\mathbf{x}))}{N}.
\end{equation*}
We have that $%
\bar H\in \overline{\mathcal{H}}_N(\Omega)$ and $\bar H \geq H$, and therefore $L_{\bar
H}\leq L_H$. On the other hand, we have for all
$x\in \Omega$,
\[
L_{\bar H}(x,u_1(x),...,u_{N-1}(x))=L_{\bar H}(x,u_1(x),...,u_{N-1}(x))+{\bar H}(x,x,...,x)\geq \sum_{\ell=1}^{N-1}\langle u_\ell(x), x\rangle.
\]
To prove (\ref{dualrep}), we use the appendix in \cite{G-M2} to deduce that for $i=2,..., N$, the gradients $\nabla_iH(x,x,...,x)$ actually exist for a.e. $ x$ in $\Omega$.\newline

The converse is straightforward since if (\ref{rep1}) holds, then (\ref{conv}%
) does, and since we also have (\ref{diago}), then the property that $%
(u_1,..., u_{N-1})$ is jointly $N$-monotone follows from (\ref{inequa}) and
the sub-antisymmetry of $H$.\hfill $\Box$\\

In the case of a single $N$-monotone vector field, we can obviously apply the above theorem to the $(N-1)$-tuple $(u,0,...,0)$ which is then $N$-monotone to find a $N$-sub-antisymmetric Hamiltonian $H$, which
is concave in the first variable, convex in the last $(N-1)$ variables such
that
\begin{equation}  \label{concave2}
\hbox{$(-u(x),u(x),0,...,0)=\nabla H(x,x,...,x)$ \quad for a.e.  $x\in
\Omega$.}
\end{equation}%
However, in this case we can restrict ourselves to $N$-cyclically sub-antisymmetric functions of two variables and establish the following extension of the Theorem of Krauss.

\begin{theorem} If $u$ is $N$-cyclically monotone on $\Omega$, then there exists a concave-convex  function of two variables $F$ that is  $N$-cyclically sub-antisymmetric and zero on the diagonal, such that
\begin{equation}
\hbox{$(-u(x), u(x))\in \partial F(x,x)$ for all $x\in \Omega$,}
\end{equation}
where $\partial H$ is the sub-differential of $H$ as a concave-convex
function \cite{Rock}. Moreover,
\begin{equation}
\hbox{$u(x)=\nabla_2 F(x,x)$ for a.e. $x\in \Omega$.}
\end{equation}

\end{theorem}
\noindent {\bf Proof:} Let $f (x,y)=\langle u(x), x-y\rangle$ and let $f^{1}\left( x,y\right) $ be
its convexification in $x$ for fixed $y$, that is
\begin{equation}
f^{1}\left( x,y\right) =\inf \left\{ \sum_{k=1}^{n}\lambda _{k}f\left(
x_{k},y\right) :\lambda _{k}\geq 0,\sum_{k=1}^{n}\lambda
_{k}=1,\sum_{k=1}^{n}\lambda _{k}x_{k}=x\right\} .
\end{equation}
Since $f(x,x)=0$, $f$ is linear in $y$, and $\sum_{i=1}^{N}f(x_{i},x_{i+1})\geq 0$ for
any cyclic family $x_{1},...,x_{N},x_{N+1}=x_{1}$ in $\Omega $, it is easy to show that  $f\geq f^1$ on $\Omega$, $f^{1}$ is convex in the first variable and concave with respect to
the second, $f^1(x,x)=0$ for each $x\in \Omega$, and that  $f^{1}$ is $N$-cyclically supersymmetric in the sense that for
any cyclic family $x_{1},...,x_{N},x_{N+1}=x_{1}$ in $\Omega $, we have
$
\sum_{i=1}^{N}f^{1}(x_{i},x_{i+1})\geq 0.
$

Consider now  $F(x, y)=-f^1(x,y)$ and note that $x\to F(x, y)$ is concave, $y\to F(x, y)$ is convex, $ F(x, y)\geq - f(x,y)=\langle u(x), y-x\rangle $ and $F$ is N-cyclically sub-antisymmetric.
By the antisymmetry, we have  \begin{equation}  \label{bounds2}
\langle u(x_1), x_2-x_1\rangle \leq F(x_1, x_2)\leq  \langle u(x_2), x_2-x_1\rangle,
\end{equation}
which yields that $(-u(x), u(x))\in \partial F(x,x)$ for all $x\in \Omega$.

Since $F$ is anti-symmetric and concave-convex, the possibly multivalued map $x\to \partial_2 F(x,x)$ is monotone on $\Omega$, and therefore single-valued and differentiable almost everywhere \cite{Ph}. This completes the proof.

\begin{remark} Note that we cannot expect to have a function $F$ such that $\sum\limits_{i=1}^{N} F (x_i, x_{i+1})=0$ for all cyclic families $x_1,..., x_N, x_{N+1}=x_1$ in $\Omega$. Actually, we believe that the only function satisfying such an $N$-antisymmetry for $N\geq 3$ must be of the form $F(x, y)=f(x)-f(y)$. This is the reason why one needs to consider functions of $N$-variables in order to get $N$-antisymmetry. In other words, the function defined by
\begin{equation}
H(x_1,x_2,..., x_N): =\frac{(N-1)F(x_1,x_2)- \sum_{i=2}^{N-1}F(x_i, x_{i+1})}{N},
\end{equation}
is $N$-antisymmetric in the sense of (\ref{exact}) and $H(x_1,x_2..., x_N) \geq F(x_1,x_2)$ for all $(x_1,x_2..., x_N)$ in $\Omega^N$.
\end{remark}

\section{Variational characterization of monotone vector fields}

In order to simplify the exposition, we shall always assume in the sequel
that $d\mu$ is Lebesgue measure $dx$ normalized to be a probability on $%
\Omega$.
 We shall also assume that $\Omega$ is convex and that its boundary has
measure zero.

\begin{theorem}
\label{var} Let $u_1,..., u_{N-1}:\Omega \to {\mathbb{R}}^d$ be bounded
measurable vector fields. The following properties are then equivalent:

\begin{enumerate}
\item The $(N-1)$-tuple $(u_1,..., u_{N-1})$ is jointly $N$-monotone a.e.,
that is there exists a measure zero set $\Omega_0$ such that $(u_1,...,
u_{N-1})$ is jointly $N$-monotone on $\Omega \setminus \Omega_0$.

\item The infimum of the following Monge-Kantorovich problem
\begin{equation}
\inf \left\{\int_{\Omega^N}\sum\limits_{\ell=1}^{N-1}\langle u_\ell(x_1),
x_1-x_{\ell +1}\rangle d\pi(x_1, x_2,..., x_N));\, \pi \in {\mathcal{P}}_{%
\mathrm{sym}}^\mu(\Omega^N) \right\}
\end{equation}%
is equal to zero, and is therefore attained by the push-forward of $\mu$ by
the map $x\rightarrow (x,x,...,x)$.

\item $(u_1,..., u_{N-1})$ is in the polar of ${\mathcal{S}}_N(\Omega, \mu)$
in the following sense,
\begin{equation}
\inf\left\{\int_\Omega \sum\limits_{\ell=1}^{N-1}\langle u_\ell (x),
x-S^\ell x\rangle \, d\mu; S \in {\mathcal{S}}_N(\Omega, \mu)\right\}=0.
\end{equation}

\item The following holds:
\begin{equation}
\inf\left\{\int_\Omega \sum\limits_{\ell=1}^{N-1}|u_\ell (x)-S^\ell x|^2
d\mu; S\in {\mathcal{S}}_N(\Omega,
\mu)\right\}=\sum\limits_{\ell=1}^{N-1}\int_\Omega |u_\ell (x)-x|^2 d\mu.
\end{equation}

\item There exists a $N$-sub-antisymmetric Hamiltonian $H$ which
is concave in the first variable, convex in the last $(N-1)$ variables, and vanishing on the diagonal such
that
\begin{equation}
\hbox{$(u_1(x),...,u_{N-1}(x))=\nabla_{2,..., N}H(x,x,...,x)$ \quad for a.e.
$x\in \Omega$.}
\end{equation}%
Moreover, $H$ is $N$-symmetric in the sense of (\ref{symmH1}).
\item The following duality holds:
\begin{equation*}
\inf \{\int_{\Omega }L_{H}(x, u_1(x),...,u_{N-1}(x))d\mu;\,H\in {\mathcal{H}}%
_{N}(\Omega )\}=\sup \{\int_{\Omega }\sum\limits_{\ell=1}^{N-1}\langle
u_\ell (x), S^\ell x\rangle \,d\mu;S\in {\mathcal{S}}_{N}(\Omega, \mu )\}
\end{equation*}%
and the latter is attained at the identity map.
\end{enumerate}
\end{theorem}

We start with the following lemma, which identifies those probabilities in ${%
\mathcal{P}}^{\mu}_{\mathrm{sym}}(\Omega^N)$ that are carried by graphs of
functions from $\Omega$ to $\Omega^N$.

\begin{lemma}
Let $S:\Omega \to \Omega$ be a $\mu$-measurable map, then the following
properties are equivalent:

\begin{enumerate}
\item The image of $\mu $ by the map $x\to (x, Sx,..., S^{N-1}x)$ belongs to
${\mathcal{P}}^{\mu}_{\mathrm{sym}}(\Omega^N)$.

\item $S$ is $\mu$-measure preserving and $S^N(x)=x$ $\mu$-a.e.

\item For any bounded Borel measurable $N$-antisymmetric $H$ on $\Omega^N$, we have $\int_\Omega H(x, Sx,...,
S^{N-1}x)\, d\mu =0$.
\end{enumerate}
\end{lemma}

\textbf{Proof.} It is clear that 1) implies 3) since $\int_{\Omega^N}H(%
\mathbf{x})\, d\pi (\mathbf{x})=0$ for any $N$-antisymmetric
Hamiltonian $H$ and any $\pi \in {\mathcal{P}}^\mu_{\mathrm{sym}}(\Omega^N)$.

That 2) implies 1) is also straightforward since if $\pi$ is the
push-forward of $\mu $ by a map of the form $x\to (x, Sx,..., S^{N-1}x)$,
where $S$ is a $\mu$-measure preserving $S$ with $S^Nx=x$ $\mu$ a.e. on $%
\Omega$, then for all $h\in L^1(\Omega^N, d\pi)$, we have
\begin{eqnarray*}
\int_{\Omega^N} h(x_1,..., x_N) d\pi&=&\int_{\Omega} h(x, Sx,...,
S^{N-1}x)\, d\mu (x) =\int_{\Omega} h(Sx, S^2x,..., S^{N-1}x, S^Nx)\, d\mu
(x) \\
&=&\int_{\Omega} h(Sx, S^2x,..., S^{N-1}x, x)\, d\mu (x) =\int_{\Omega^N}
h(\sigma (x_1,..., x_N)) d\pi.
\end{eqnarray*}
We now prove that 2) and 3) are equivalent. Assuming first that $S$ is $\mu$%
-measure preserving such that $S^N=I$ $\mu$ a.e., then for every Borel
bounded $N$-antisymmetric $H$, we have
\begin{eqnarray*}
\int_\Omega H(x, Sx, S^2x,...,S^{N-1}x)d\mu&=&\int_\Omega H(Sx,
S^2x,...,S^{N-1}x, x)d\mu \\
&=&...=\int_\Omega H(S^{N-1}x, x, Sx,...,S^{N-2}x)d\mu.
\end{eqnarray*}
Since $H$ is $N$-antisymmetric, we can see that
\begin{equation*}
H(x, Sx,...,S^{N-1}x)+H(Sx, S^2x,...,S^{N-1}x, x)+... H(S^{N-1}x, x, Sx,..,
S^{N-2}x)=0.
\end{equation*}
It follows that $N\int_\Omega H(x, Sx, S^2x,.., S^{N-1}x)d\mu=0$.\newline

For the reverse implication, assume $\int_\Omega H (x, S x,S^2x, ...,
S^{N-1}x) d\mu=0$ for every $N$-antisymmetric Hamiltonian $H$. By
testing this identity with the Hamiltonians
\begin{equation*}
H (x_1, x_2,..., x_N)=f(x_1)-f(x_i),
\end{equation*}
where $f$ is any continuous function on $\Omega$, one gets that $S$ is $\mu$%
-measure preserving. Now take the Hamiltonian
\begin{equation*}
H (x_1, x_2,..., x_N)=|x_1-Sx_N|-|S x_1-x_2|-|x_{2}-S x_1|+|Sx_2-x_{3}|.
\end{equation*}
Note that $H \in {\mathcal{H}}_N(\Omega)$ since it is of the form $H
(x_1,..., x_N)=f(x_1,x_2, x_N)-f(x_2, x_3, x_1)$. Now test the above
identity with such an $H $ to obtain
\begin{eqnarray*}
0=\int_\Omega H (x, S x,S^2x, ..., S^{N-1}x) d\mu =\int_\Omega
|x-SS^{N-1}x|\, d\mu.
\end{eqnarray*}
It follows that $S^N=I$ $\mu$ a.e. on $\omega$, and we are done.\hfill $\Box$%
\bigskip

\noindent\textbf{Proof of Theorem \ref{var}:} To show that (1) implies (2),
it suffices to notice that if $\pi $ is a $\sigma$-invariant probability
measure on $\Omega^N$ such that $\mathrm{proj}_{1} \pi =\mu$, then
\begin{eqnarray*}
\int_{\Omega ^{N}}\sum\limits_{\ell=1}^{N-1}\langle u_\ell (x_1),
x_1-x_{\ell +1}\rangle d\pi \left( x_{1},...,x_N\right)&=&\frac{1}{N}%
\sum_{i=1}^N \int_{\Omega ^{N}}\sum\limits_{\ell=1}^{N-1}\langle
u_\ell\left( x_{i}\right) ,x_{i}-x_{i+\ell}\rangle d\pi \left(
x_{1},...,x_N\right) \\
&=&\frac{1}{N}\int_{\Omega ^{N}}
\left(\sum_{i=1}^N\sum\limits_{\ell=1}^{N-1} \left\langle u_\ell \left(
x_{i}\right) ,x_{i}-x_{i+\ell}\right\rangle\right) d\pi \left(
x_{1},...,x_N\right) \\
&\geq&0,
\end{eqnarray*}
since $(u_1,..., u_{N-1})$ is jointly $N$-monotone. On the other hand, if $%
\pi$ is the $\sigma$-invariant measure obtained by taking the image of $%
\mu:=dx$ by $x\to (x,..., x)$, then
\begin{equation*}
\int_{\Omega ^{N}}\sum\limits_{\ell=1}^{N-1}\langle u_\ell (x_1),
x_1-x_{\ell +1}\rangle d\pi \left( x_{1},...,x_N\right)=0.
\end{equation*}
To show that (2) implies (3), let $S$ be a $\mu$-measure preserving
transformation on $\Omega$ such that $S^{N}=I$ $\mu$ a.e. on $\Omega$. Then
the image $\pi _{S}$ of $\mu$ by the map%
\begin{equation*}
x\rightarrow \left( x,Sx,S^{2}x,..., S^{N-1}x\right)
\end{equation*}%
is $\sigma$-invariant, hence
\begin{equation*}
\int_{\Omega^N}\sum\limits_{\ell=1}^{N-1}\langle u_\ell (x_1), x_1-x_{\ell
+1}\rangle d\pi_S \left( x_{1},...,x_N\right)=\int_\Omega
\sum\limits_{\ell=1}^{N-1}\langle u_\ell (x), x-S^\ell x\rangle \, d\mu\geq
0.
\end{equation*}
By taking $S=I$, we get that the infimum is necessarily zero.\newline

The equivalence of (3) and (4) follows immediately from developing the
square.
\newline

We now show that (3) implies (1). For that take $N$ points $%
x_{1},x_{2},...,x_{N}$ in $\Omega$, and let $R>0$ be such that $B\left(
x_{i},R\right) \subset \Omega $. Consider the transformation
\begin{equation*}
S_{R}\left( x\right) =\left\{
\begin{array}{c}
x-x_{1}+x_{2}\text{ for }x\in B\left( x_{1},R\right) \\
x-x_{2}+x_{3}\text{ for }x\in B\left( x_{2},R\right) \\
... \\
x-x_{N}+x_{1}\text{ for }x\in B\left( x_{N},R\right) \\
x\text{ otherwise}%
\end{array}%
\right.
\end{equation*}
It is easy to see that $S_R$ is a measure preserving transformation and that
$S_R^{N}=Id$. We then have
\begin{eqnarray*}
0 \leq \int_{\Omega }\sum\limits_{\ell=1}^{N-1}\langle u_\ell (x),
x-S_R^\ell x\rangle\, d\mu \leq \sum_{i=1}^{N}\int_{B\left( x_{i},R\right)
}\sum\limits_{\ell=1}^{N-1}\left\langle u_\ell\left( x \right)
,x_{i}-x_{\ell +i}\right\rangle d\mu.
\end{eqnarray*}
Letting $R\rightarrow 0$, we get from Lebesgue's density theorem, that
\begin{equation*}
\frac{1}{\left\vert B\left( x_{i},R\right) \right\vert }\int_{B\left(
x_{i},R\right) }\left\langle u_\ell \left( x\right)
,x_{i}-x_{\ell+i}\right\rangle d\mu\rightarrow \left\langle u_\ell \left(
x_{i}\right) ,x_{i}-x_{\ell+i}\right\rangle,
\end{equation*}%
from which follows that $(u_1,..., u_{N-1})$ are jointly $N$-monotone a.e. on $\Omega$.
The fact that (1) is equivalent to (5) follows immediately from Theorem \ref%
{main.1}.

\noindent To prove that 5) implies 6) note that for all $p_i\in {\mathbb{R}}%
^d, x\in \Omega, y_i\in \Omega, i=1,...,N-1$,
\begin{equation*}
L_H(x,p_1,..., p_{N-1})+H(x, y_1,..., y_{N-1}) \geq \sum_{i=1}^{N-1}\langle
p_i, y_i\rangle,
\end{equation*}
which yields that for any $S\in {\mathcal{S}}_N(\Omega, \mu)$,
\begin{equation*}
\int_\Omega [ L_H(x,u_1(x),...,u_{N-1}(x))\, d\mu + H(x, Sx,...,
S^{N-1}x)]\, d\mu \geq \int_\Omega \sum\limits_{\ell=1}^{N-1}\langle u_\ell
(x), S^\ell x\rangle \, d\mu.
\end{equation*}
If $H \in {\mathcal{H}}_N(\Omega)$ and $S\in {\mathcal{S}}_N(\Omega, \mu)$,
we then have $\int_\Omega H(x, Sx,..., S^{N-1}x) d\mu =0$, and therefore
\begin{equation*}
\int_\Omega L_H(x,u_1(x),...,u_{N-1}(x))\, d\mu \geq \int_\Omega
\sum\limits_{\ell=1}^{N-1}\langle u_\ell (x), S^\ell x\rangle \, d\mu.
\end{equation*}
If now $H$ is the $N$-sub-antisymmetric Hamiltonian obtained by 5), which is
concave in the first variable, convex in the last $(N-1)$ variables,
then
\begin{equation*}
\hbox{$ L_H(x,u_1(x),...,u_{N-1}(x)) +
H(x,x,...,x)=\sum_{\ell=1}^{N-1}\langle u_\ell(x), x\rangle$ \, for all
$x\in \Omega\setminus \Omega_0$, }
\end{equation*}
and therefore $\int_\Omega L_H(x,u_1(x),...,u_{N-1}(x)) \,
d\mu=\sum_{\ell=1}^{N-1}\int_\Omega \langle u_\ell(x), x\rangle\, d\mu$.

Consider now
\begin{equation*}
\bar H (\mathbf{x})=\frac{(N-1)H(\mathbf{x})- \sum_{i=1}^{N-1}H(\sigma^i(%
\mathbf{x}))}{N}.
\end{equation*}
As before,
we have that $%
\bar H\in {\mathcal{H}}_N(\Omega)$ and $\bar H \geq H$. Since $L_{\bar
H}\leq L_H$, we have that $\int_\Omega L_{\bar H}(x,u_1(x),...,u_{N-1}(x))\, d\mu
= \sum_{\ell=1}^{N-1}\int_\Omega \langle u_\ell(x), x\rangle\, d\mu$ and (6)
is proved.\newline

Finally, note that (6) readily implies (3), which means that $(u_1,...,
u_{N-1})$ is then jointly $N$-monotone. \hfill $\Box$\\

We now consider again the case of a single $N$-cyclically monotone vector field.

\begin{corollary}
Let $u:\Omega \to {\mathbb{R}}^d$ be a bounded measurable vector field. The
following properties are then equivalent:

\begin{enumerate}
\item $u$ is $N$-cyclically monotone a.e., that is there exists a measure
zero set $\Omega_0$ such that $u$ is $N$-cyclically monotone on $\Omega
\setminus \Omega_0$.

\item The infimum of the following Monge-Kantorovich problem
\begin{equation}
\hbox{$\inf\{\int_{\Omega^N}\langle u(x_1), x_1-x_2\rangle d\pi({\bf x});\,
\pi \in {\mathcal P}_{\rm sym}^\mu(\Omega^N)$}\}
\end{equation}%
is equal to zero, and is therefore attained by the push-forward of $\mu$ by
the map $x\rightarrow (x,x,...,x)$.

\item The vector field $u$ is in the polar of ${\mathcal{S}}_N(\Omega, \mu)$%
, that is
\begin{equation}
\inf\{\int_\Omega \langle u(x), x-Sx\rangle \, d\mu; S \in {\mathcal{S}}%
_N(\Omega, \mu)\}=0.
\end{equation}

\item The projection of $u$ on ${\mathcal{S}}_N(\Omega, \mu)$ is the
identity map, that is
\begin{equation}
\inf\{\int_\Omega |u(x)-Sx|^2 d\mu; S\in {\mathcal{S}}_N(\Omega,
\mu)\}=\int_\Omega |u(x)-x|^2 d\mu.
\end{equation}

\item There exists a $N$-cyclically sub-antisymmetric function  $H$ of two variables, which is concave in the first variable, convex in the second  variable, vanishing on the diagonal and
 such
that
\begin{equation}  \label{concave}
\hbox{$u(x)=\nabla_2 H(x,x)$ \quad for a.e. $x\in
\Omega$.}
\end{equation}%

\item The following duality holds:
\begin{equation*}
\inf \{\int_{\Omega }L_{H}(x,u(x),0, ...,0)d\mu;\,H\in {\mathcal{H}}%
_{N}(\Omega )\}=\sup \{\int_{\Omega }\langle u(x),Sx\rangle \,d\mu;S\in {%
\mathcal{S}}_{N}(\Omega, \mu )\}
\end{equation*}%
and the latter is attained at the identity map.
\end{enumerate}
\end{corollary}

\noindent\textbf{Proof:} This is an immediate application of Theorem \ref%
{var} applied to the $(N-1)$-tuplet vector fields $(u,0,..., 0)$, which is
clearly jointly $N$-monotone on $\Omega \setminus \Omega_0$, whenever $u$ is
$N$-monotone on $\Omega \setminus \Omega_0$.

\begin{remark} \rm
\textrm{Note that the sets of $\mu$-measure preserving $N$-involutions $({%
\mathcal{S}}_N(\Omega, \mu))_N$ do not form a nested family, that is ${%
\mathcal{S}}_N(\Omega, \mu)$ is not necessarily included in ${\mathcal{S}}%
_M(\Omega, \mu)$, whenever $N\leq M$, unless of course $M$ is a multiple of $%
N$. On the other hand, the above theorem shows that their polar sets, i.e.,
\begin{equation*}
\hbox{${\mathcal S}_N(\Omega, \mu)^0=\{u\in L^2(\Omega, \R^d);\, \int_\Omega
\langle u(x), x-Sx\rangle \, d\mu \geq 0$ for all  $S \in
{\mathcal{S}}_N(\Omega, \mu)\}$},
\end{equation*}
which coincide with the $N$-cyclically monotone maps, satisfy
\begin{equation*}
{\mathcal{S}}_{N+1}(\Omega,\mu)^0 \subset {\mathcal{S}}_N(\Omega, \mu)^0,
\end{equation*}
for every $N\geq 1$. This can also be seen directly. Indeed, it is clear
that a $2$-involution is a $4$-involution but not necessarily a $3$%
-involution. On the other hand, assume that $u$ is $3$-cyclically monotone
operator, then for any transformation $S:\Omega \to \Omega$, we have
\begin{equation*}
\int_\Omega \langle u(x), x-Sx\rangle d\mu+\int_\Omega \langle u(Sx),
Sx-S^2x\rangle d\mu+\int_\Omega \langle u(S^2x), S^2x-x\rangle d\mu\geq 0.
\end{equation*}
If now $S$ is measure preserving, we have
\begin{equation*}
\int_\Omega \langle u(x), x-Sx\rangle d\mu+\int_\Omega \langle u(x),
x-Sx\rangle d\mu+\int_\Omega \langle u(S^2x), S^2x-x\rangle d\mu\geq 0,
\end{equation*}
and if $S^2=I$, then $\int_\Omega \langle u(x), x-Sx\rangle d\mu \geq 0$,
which means that $u \in {\mathcal{S}}_2(\Omega,\mu)^0$. Similarly, one can
show that any $(N+1)$-cyclically monotone operator belongs to ${\mathcal{S}}%
_N(\Omega, \mu)^0$. In other words, ${\mathcal{S}}_{N+1}(\Omega,
\mu)^0\subset {\mathcal{S}}_N(\Omega, \mu)^0$ for all $N\geq 2$. Note that ${%
\mathcal{S}}_1(\Omega, \mu)^0=\{I\}^0=L^2(\Omega, {\mathbb{R}}^d)$, while
\begin{equation*}
\hbox{${\mathcal S}(\Omega, \mu)^0= \cap_N {\mathcal S}_N(\Omega, \mu
)^0=\{u\in L^2(\Omega, \R^d), u=\nabla \phi$ for some convex function $\phi$
in $W^{1,2}(\R^d)$}\},
\end{equation*}
in view of classical results of Rockafellar \cite{Rock} and Brenier \cite{Br}%
. }
\end{remark}

\begin{remark} \rm
In a forthcoming paper \cite{G-M2}, the above result is extended to give a
similar decomposition for any family of bounded measurable vector fields $%
u_{1},u_{2},....,u_{N-1}$ on $\Omega $. It is shown there that there exists
a measure preserving $N$-involution $S$ on $\Omega $ and an $N$-antisymmetric Hamiltonian $H$ on $\Omega ^{N}$ such that for $i=1,...,N-1$,
we have
\begin{equation*}
\hbox{$u_i(x)=\nabla_{i+1}H(x, Sx, S^2x,...S^{N-1}x)$\quad  for a.e. $x\in
\Omega.$}
\end{equation*}%
\end{remark}
{\bf Acknowledgement:} We are grateful to the anonymous referee for a careful reading of this paper, which led to several improvements.


\begin{thebibliography}{99}
\bibitem{Br} Y. Brenier, \emph{Polar factorization and monotone
rearrangement of vector-valued functions,} Comm. Pure Appl. Math. \textbf{44}
(1991), 375-417.

\bibitem{13} S. P. Fitzpatrick, \textit{Representing monotone operators by
convex functions}, Proc. Centre for Math. Analysis \textbf{20} (1989), 59-65.

\bibitem{Gg} W. Gangbo, \emph{An elementay proof of the polar factorization
of vector-valued functions}, Arch. Rat. Math. Analysis \textbf{128}, No.5,
(1994) 381-399.

\bibitem{Gh} N. Ghoussoub, \textit{Selfdual partial differential systems and
their variational principles,} Springer Monograph in Mathematics,
Springer-Verlag (2008), 356 p.


\bibitem{Gh-Ma} N. Ghoussoub, B. Maurey, \textit{Remarks on multidimensional
symmetric Monge-Kantorovich problems}, Discrete and Continuous Dynamical
Systems-A, To appear (2013)

\bibitem{G-M2} N. Ghoussoub, A. Moameni, \textit{Symmetric Monge-Kantorovich
problems and polar decompositions of vector fields}, Preprint (2012)


\bibitem{G-M} N. Ghoussoub, A. Moameni, \textit{A Self-dual Polar
Factorization for Vector Fields}, Comm. Pure. Applied. Math., Vol 66, Issue 6 (2013) p. 905-933

\bibitem{Mi} P. Millien, \textit{A functional analytic approach to the
selfdual polar decomposition}, Master thesis, UBC (2011).

\bibitem{Kra} E. Krauss, \emph{A representation of arbitrary maximal
monotone operators via subgradients of skew-symmetric saddle functions,}
Nonlinear Anal. 9 (1985), no. 12, 1381-1399,

\bibitem{Ph} R. R. Phelps, \textit{Convex functions, monotone operators and
differentiability}, Lecture Notes in Math. 1364, Springer Verlag, New York,
Berlin, Tokyo, (1998), 2nd edition 1993.

\bibitem{Rock} T. Rockafellar, \emph{Convex Analysis,} 1970, Princeton
University Press.
\end{thebibliography}
\end{document}